\documentclass[12pt]{article}

\usepackage{graphicx}
\usepackage{amssymb}
\usepackage{amsmath}
\usepackage{color}

\newcommand{\kom}[1]{}
\renewcommand{\kom}[1]{{\bf [#1]}}
\definecolor{blau}{rgb}{0.1,0.0,0.9}

\newcounter{komcounter}
\numberwithin{komcounter}{section}

\begin{document}

\title{\bf How to find simple nonlocal stability and resilience measures}
\author{Niklas L.P. Lundstr\"{o}m
 \linebreak \\\\
\it \small Department of Mathematics and Mathematical Statistics, Ume{\aa} University\\
\it \small SE-90187 Ume{\aa}, Sweden\/{\rm ;}
\it \small niklas.lundstrom@umu.se\linebreak \\\\}

\maketitle

\begin{abstract}
Stability of dynamical systems is a central topic with applications in widespread areas such as
economy, biology, physics and mechanical engineering.
The dynamics of nonlinear systems may completely change due to perturbations forcing the solution to jump from a safe state into another,
possibly dangerous, attractor.
Such phenomena can not be traced by the widespread local stability and resilience measures,
based on linearizations, 
accounting only for arbitrary small perturbations.
Using numerical estimates of the size and shape of the basin of attraction,
as well as the systems returntime to the attractor after given a perturbation,
we construct simple nonlocal stability and resilience measures that
record a systems ability to tackle both large and small perturbations.
We demonstrate our approach
on the Solow-Swan model of economic growth,
an electro-mechanical system as well as on
a stage-structured population model,
and conclude that the suggested measures detect dynamic behaviour,
crucial for a systems stability and resilience,
which can be completely missed by local measures.
The presented measures are also easy to implement on a standard laptop computer.
We believe that our approach 
will constitute an important step towards filling a current gap in the literature by putting forward and explaining simple ideas and methods,
and by delivering explicit constructions of several promising nonlocal stability and resilience measures.\\
\end{abstract}



\setcounter{equation}{0} \setcounter{theorem}{0}

\section{Introduction}
\label{sec:intro}


Understanding stability of dynamical systems is important for widespread areas of research such as
e.g.~economy, biology, physics and mechanical engineering. 
As computers and computational tools have become more and more efficient,
simulations and numerical investigations of realistic, high-dimensional, advanced mathematical models have become easier
and also increasingly popular.
As such advanced models account for numerous factors and their interplay,
the dynamics are often nontrivial, nonlinear and coexisting attractors (bistability) exists or may be hard to rule out.
Such systems algebraic structure are also often complicated
and classical mathematical analysis is often difficult.
This calls for further research in mathematical analysis,
as well as in numerical methods,
for finding efficient ways of quantifying the stability of dynamical systems.

\begin{figure}[h]
\begin{center}
\includegraphics[height = 4cm, width = 13cm]{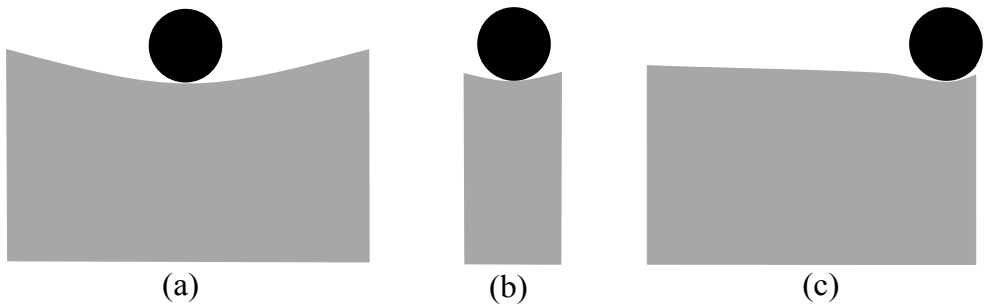}
\caption{Three identical balls rest at equilibrium on three different stands.
All stands have the same shape in a neighborhood of the equilibrium and hence a
local approach would rank the systems in (a), (b) and (c) as equally stable and equally resilient.
The basin of attraction is largest in (a), followed by (c) while it is smallest in (b).
Therefore, measuring only the size of the basin of attraction would rank (b) as least stable
even though it is easier to push the ball off the stand in (c) as the equilibrium in (c) is closer to the boundary of the basin.
Hence, both the size and the shape of the basin of attraction should be included when measuring stability and resilience of a nonlinear dynamical system.}
\label{fig:fig_1}
\end{center}
\end{figure}

The analysis of the stability of attractors (e.g.~equilibrium points, limit cycles, quasiperiodic or strange attractors)
in dynamical systems naturally split into local analysis and nonlocal analysis.
The local stability approach is usually based on linearizations and yields information in a small neighbourhood of the attractor,
saying little or nothing about the systems behaviour a bit away from the attractor. 
Measuring stability in nonlinear dynamical systems using local methods
(such as eigenvalues of the Jacobian matrix at an equilibrium)
delivers information of how the system reacts only on arbitrary small perturbations.
A perturbation of a given size may push a locally stable state into another attractor having a completely different,
possibly dangerous, behaviour.
For example, researchers believe that the crash of the aircraft YF-22 Boeing in April 1992 was caused by a sudden switch to an unsafe attractor (Dudkowski et al 2016; Lauvdal et al 1997).
Figure \ref{fig:fig_1} shows three systems having completely different abilities to withstand perturbations,
but exactly the same local stability and local resilience.

The local stability analysis should therefore be complemented by, or replaced by,
a nonlocal approach that considers properties of the basin of attraction
(henceforth \textit{basin})
to the attractor under investigation.
The \emph{size of basin} (it's volume) constitutes a natural candidate for a nonlocal stability measure as a large basin indicates that the system comes back to
the attractor with a high probability, given a random perturbation. It is also easy to numerically estimate the size of the basin.
However,
if the distance from the attractor to the boundary of the basin is short in some direction,
then a small perturbation in this direction can push the system to another,
possibly dangerous, attractor, even though the basin is large,
see Figure \ref{fig:fig_1} (c).
Therefore, the \emph{shape of basin} is another natural candidate to be included in a nonlocal measure construction.
It is not trivial how to construct measures based on the shape of the basin,
but we may conclude that a large and convex basin,
having the attractor in the middle,
should indicate high stability of the system.
This is because in such case a large perturbation can be given to the system,
in any direction,
without pushing the system to another attractor.

In the field of mathematics, it is classical to attempt to find Lyapunov functions in order to prove nonlocal stability of attractors in dynamical systems.
Such method gives analytical estimates of both the size and the shape of the basin.
However, it is often a difficult task to find suitable Lyapunov functions,
making this approach insufficient,
especially when dealing with high-dimensional complicated systems of equations.
If a Lyapunov function is not available, another investigation of the basin has to be considered in order to
understand the nonlocal stability of the attractor.
The aim of this paper is to deliver simple methods for doing this by considering several nonlocal stability and resilience measures built upon
numerical estimates of the basin and on the returntime of trajectories corresponding to perturbations from the attractor.

Keeping in mind that works involving nonlocal stability approaches can be hard to find
as papers generally put focus on the results of a stability analysis rather than the stability measure itself,
a literature survey 
indicates that nonlocal stability measures are rarely used. 
Most works connecting to applied dynamical systems,
e.g.~in theoretical ecology and mechanical engineering,
seem to rely solely on a local investigation of stability.
This is surprising because even though a mathematical approach
is difficult, today's computers can often calculate estimates of the basin within a couple of minutes.
Lundstr\"om (2006) and Lundstr\"om and Aidanp\"a\"a 
(2007) used two nonlocal stability measures 
when investigating stability of an electric generator. 
The authors refer to the systems ability to absorb perturbations as \emph{robustness} of the system.
Their first stability measure estimates the basin's size,
%
%
reflecting the probability that the generator system returns to the ``safe" attractor given a perturbation.
A similar approach was used by Menck et al.
(2013) 
to prove properties of neutral networks and power grids;
they refer to the nonlocal stability measure as \emph{basin stability}.
Menck et al. (2013) also give several examples in which local stability measures fail to detect crucial destabilizing effects.
Lundstr\"om and Aidanp\"a\"a (2007) used also a stability measure based upon the shape of the basin through approximation of the smallest perturbation needed to push the system from the safe attractor into another ``unsafe" attractor.
A similar idea was suggested by Klinshov et al. (2016) 
and algorithms for calculating such measure in the setting of general attractors exists
(Kerswell et al. 2014; Klinshov et al. 2016).

A first aim of this paper is to expand on the above ideas 
by constructing and evaluating simple nonlocal measures of stability using fundamental properties of the basin such as it's size and shape.
To do so we 
consider two nonlocal measures accounting for both large and small perturbations.
The first measure, denoted by $\mathcal{P}$, estimates the size of the basin and thereby answers the natural question; \emph{what is the probability that the system returns to the attractor given a perturbation?}
The second measure, denoted by $\mathcal{D}$, builds on the shape of the basin and concerns the question; \emph{what is the smallest perturbation needed to push the solution to another attractor?}
The measure $\mathcal{D}$ calls for suitable ways to compare distances in the state-space, and, focusing on mechanical systems,
we construct the version $\mathcal{D}_{\text{energy}}$ reflecting the least amount of energy needed to push the systems solution into another attractor.
We also suggest the version $\mathcal{D}_{\text{rel}}$, considering relative distances,
and show how to use a version of it when evaluating harvesting strategies in stage-structured populations.

A second aim of this paper is to include a next natural candidate for nonlocal stability into the measures construction,
i.e.~the \emph{returntime} of a trajectory to the attractor given a perturbation.
Doing so leads us into the \emph{resilience} of a system and construction of resilience measures;
an important concept increasingly used in ecology
(see e.g.~Neubert and Caswell 1997; Loreau and Behera 1999; Loeuille 2010; Arnoldi et al. 2016; Haegeman et al. 2016; Arnoldi et al. 2017).
There are at least two different definitions of resilience of a system in the literature.
The first definition, and the more traditional, concentrates on stability near an equilibrium,
where resistance to disturbance and speed of return to the equilibrium, following a perturbation, are used to measure the property
(O'Neill et al., 1986; 
Pimm 1984; 
Tilman and Downing 1994). 
This view of resilience provides one of the foundations for economic theory as well and may be termed engineering resilience (Holling 1996). 
The second definition emphasizes conditions far from any equilibrium,
where instabilities can flip a system into another regime of behavior,
i.e., to another stability domain (Holling 1973). 
In this case the measurement of resilience is the magnitude of disturbance that can be absorbed before the system changes its structure by changing the variables and processes that control behavior.
This second view may be termed ecological resilience (Walker et al. 1969) and clearly involves properties of the basin of attraction.

Even though neither of the above definitions of resilience restrict to small perturbations,
widespread technics for calculating resilience are
based solely on local analysis by considering e.g.~the eigenvalues of the Jacobian matrix at equilibrium.
See e.g.~Neubert and Caswell (1997) 
and Arnoldi et al. (2016) 
and the references therein for more on standard local resilience measures, 
alternative resilience measures and their relations as well as
discussions concerning shortcomings of the local approach in general.
For further discussion on the use of local resilience in ecology and the fact that it can be difficult to assess from an
empirical point of view, see
Haegeman et al. (2016).
Mitra et al. (2015) 
suggest an integrative measure of resilience based on both local and nonlocal features,
and Lundstr\"om et al. (2016) 
use a simple nonlocal resilience measure for evaluating harvesting strategies 
of age- and stage-structured populations.
However, our impression is that,
in analogue with the literature on pure stability measures discussed above,
the nonlocal approach is under-used also in the setting of resilience.
As a consequence, valuable information of a systems ability to sustain perturbations may be unrecorded:
Recall that all three systems in Figure \ref{fig:fig_1} have identical resilience according to any measure based on a
local linearization (such as the largest eigenvalue) near the equilibrium.

This motivates our second aim;
to find simple nonlocal resilience measures based on the attractors basin as well as on the
returntime for trajectories starting within the basin
(and thus returns to the attractor under investigation).
We will present four nonlocal resilience measures, accounting for large as well as small perturbations,
of which the first two steams from seeing resilience as the reciprocal of the returntime, i.e. as a rate of return given a perturbation.
Indeed, our simple measure $\mathcal{R}$ will answer the question; \emph{what is the expected rate of return of the system given a random perturbation?}
We also keep track of a ``worst-case resilience", or slowest recovery from a set of perturbations,
through the measure
$\mathcal{R}_{\text{worst}}$ answering; \emph{what is the slowest rate of return of the system given a random perturbations?}

We proceed by introducing the concept \emph{basin-time} as the subset of the basin of attraction from which all trajectories return to the attractor within a certain time-limit.
Based upon the basin-time we build two more resilience measures.
The first is $\mathcal{P}^{\tau}$ which concerns the question; \emph{what is the probability that the system returns from a random perturbation in $\tau$ years?}
The second is $\mathcal{D}^{\tau}$; \emph{what is the least perturbation from which the system will not return in $\tau$ years?}
The above questions should be relevant for wide ranges of applications e.g.~in biology, fisheries management, economy, and management of financial systems.
Our resilience measures are applicable in the setting of simple linear systems as well as in high-dimensional nonlinear systems.
All suggested stability and resilience measures are easy to implement on a standard laptop computer, and, together,
they yield an explaining picture of the systems dynamics in a ``large" neighborhood of the attractor.

In Section \ref{sec:methods} we present the main ideas and definitions of the nonlocal stability and resilience measures,
while in Section \ref{sec:results} we demonstrate
our approach on three independent mathematical models;
the Solow-Swan model of economic growth, 
a simple electro-mechanical model and 
a stage-structured population model.
In each application we discuss what can be learned by each measure, and why.
We also consider relations between the measures as well as their relation to linear measures.
In Section \ref{sec:discussion} we discuss the achieved results,
how to apply the simple measures efficiently to more general situations, 
some topics for future research and finally we give some concluding words.


\section{Methods}
\label{sec:methods}
%
We focus on models governed by a dynamical system 
on the standard form
\begin{align}\label{eq:dyn-syst}
\frac{\partial x}{\partial t} \,=\, f(x,t), \quad  x \in \mathbf{R}^n,\; t > 0 \quad \text{and initial condition}\quad x(0) = \bar x,
\end{align}
where $f(x,t)$ is a given function defined on $\mathbf{R}^n\times\mathbf{R}$ taking values in $\mathbf{R}^n$ and $\bar x \in \mathbf{R}^n$. 
%
However, our ideas are general and the following measure constructions have potential for generalization to systems described by a wider range of mathematical tools as well,
such as e.g.~systems of partial differential equations,
systems combining partial and ordinary differential equations,
stochastic differential equations and
cellular automaton.
%

%

To set up the measures the fist step is to choice a suitable set of perturbations to test the system for.
This choice 
may preferably reflect what the system may be exposed to in reality.
Perturbations can be taken deterministically or randomly from a predefined probability distribution;
the measures presented below can handle any reasonable choice of perturbations.
A perturbation is inferred through an initial conditions of \eqref{eq:dyn-syst}
and thus the set of perturbation is a set of initial conditions.

For each initial condition in the set of chosen perturbations,
the corresponding trajectory is numerically integrated in order to see if it returns to the attractor
(of which stability and/or resilience are to be investigated) or not.
All initial conditions that returned are saved as ``safe" in a set $\mathbb{I}_{\text{safe}}$
and those that not returned as ``unsafe" in a set $\mathbb{I}_{\text{unsafe}}$.
The time needed for the trajectory to recover the attractor is also saved as the returntime $T(x_i)$ for each safe initial condition $x_i \in \mathbb{I}_{\text{safe}}$.
Let $\mathbb{N}_{\text{safe}}$, $\mathbb{N}_{\text{unsafe}}$ and $\mathbb{N}_{\text{tot}}$ denote the safe,
the unsafe and the total number of tested initial conditions,
respectively.
By ``return to the attractor" we mean that the trajectory has entered a small pre-defined neighborhood of the attractor.
The choice of this (very) small neighborhood may be based on what is considered as normal, or perfect, operation of the system.

Even though different sets of perturbations of course can be considered for each measure,
we highlight that we are able to get estimates of all suggested measures from one single sample of perturbations.
Since most computational time lies in the integration of each trajectory,
it is therefore no more expensive to estimates all measures than just one single measure.

\subsubsection*{Nonlocal stability measures}

We choice to measure nonlocal stability through the \emph{size of basin} 
using the simple measure $\mathcal{P}$ estimated by
\begin{align*}
{\mathcal{\hat P}} = \frac{\mathbb{N}_{\text{safe}}}{\mathbb{N}_{\text{tot}}}.
\end{align*}
The measure $\mathcal{P}$ reflects the probability that the system returns to the attractor given a perturbation from the pre-specified set of initial conditions.
We have $0 \leq \mathcal{P} \leq 1$; if $\mathcal{P} = 0$ then no trajectories return and all perturbations are unsafe,
and in case all trajectories return then $\mathcal{P} = 1$.
The similar measure was used by 
Lundstr\"om and Aidanp\"a\"a (2007) for estimating the stability of electric generators and by Menck et al. (2013) for proving results on network stability.

We measure stability through the \emph{shape of basin} of an attractor by considering
the shortest distance from the attractor to the boundary of the basin.
We denote this measure by $\mathcal{D}$ and estimate it as
\begin{align*}
\mathcal{\hat D} = \min_{\substack{x \,\in \,\mathbb{I}_{\text{unsafe}}}} \text{dist}(x, \mathbb{A}),
\end{align*}
%
%
%
%
%
where $\mathbb{A}$ is the attractor of which stability is to be investigated and
$\text{dist}(x, \mathbb{A})$ denotes the distance from $x$ to $\mathbb{A}$ defined in a suitable way.
The measure $\mathcal{D}$ reveals the smallest perturbation needed to push the system to another attractor.
A similar measure was considered by Klinshov et al. (2016) 
and algorithms for calculating it in the setting of general attractors exists
(Kerswell et al. 2014; Klinshov et al. 2016).
In Section \ref{sec:results} we will demonstrate our approach to estimate $\mathcal{D}$ in the setting of equilibria,
and in Section \ref{sec:discussion} we further discuss how one can estimate our suggested measures for some general attractors in a simple way.

It is important and often nontrivial to find a suitable distance function (norm) for the specific application of the measure $\mathcal{D}$.
Standard Euclidean distance can be used due to its simplicity,
but is seldom motivated as different dimensions in the state-space may have completely different meanings in terms of perturbations,
e.g.~in mechanical systems we end up with comparing velocity-dimensions to displacement-dimensions.
To deal with such case, we consider an energy-norm to construct the version $\mathcal{D}_{\text{energy}}$,
which finds the \emph{least amount of energy} needed to push push the system from its current state into another attractor.
Indeed, a perturbation can be given through a displacement, a velocity impulse, or a combination of these.
The energy-norm is defined as the potential energy needed to perform the displacement perturbation plus the kinetic energy given through the velocity impulse. We further explain and demonstrate the idea in subsection \ref{sec:mech}.
When dealing with biomass in population dynamics we found it relevant to use a relative distance function and introduce the version $\mathcal{D}_{\text{rel}}$,
considering the size of the perturbation in relation to the biomass at the original state.
We explain and demonstrate this idea in subsection \ref{sec:biol}.
%
As the measure $\mathcal{D}$ searches for a particular perturbation, ``the most dangerous one,"
one may consider a distribution of perturbations that puts more weight near the boundary of the basin (if such information is available)
in order to estimate it more efficiently.

\subsubsection*{Nonlocal resilience measures}

By seeing resilience as the reciprocal of the returntime we consider the simple resilience measure $\mathcal{R}$ through the estimate
\begin{align*}
%
\mathcal{\hat R} = \frac{1}{\mathbb{N}_{\text{tot}}} \, \sum_{x_i \,\in \,\mathbb{I}_{\text{safe}}} \, \frac{1}{T(x_i) + t_{\epsilon}}.
\end{align*}
The measure $\mathcal{R}$ accounts for the basin of attraction by giving zero resilience to initial conditions in $\mathbb{I}_{\text{unsafe}}$
of which trajectories did not return to the attractor.
The parameter $t_{\epsilon} > 0$ is added to avoid division by zero and explosion for very small perturbations (resulting in zero returntime),
and can often be chosen to 1 (with exception of very fast systems where it may be taken smaller).
The measure $\mathcal{R}$ yields the expected rate of return of the system given a random perturbation.

We consider also the worst-case resilience, $\mathcal{R}_{\text{worst}}$,
as the reciprocal of the slowest recovery (longest returntime) through the estimate
\begin{align*}
\mathcal{\hat R}_{\text{worst}} =\, \min_{\substack{x \,\in \,\mathbb{I}_{\text{safe}}}} \frac{1}{T(x) + t_{\epsilon}}.
\end{align*}
%
%
The measure $\mathcal{R}_{\text{worst}}$ reflects the slowest rate of return of the system given a random perturbation.
It is often relevant to restrict perturbations further in order to specify this measure.
This can be done simply by replacing the set of all safe perturbations, $\mathbb{I}_{\text{safe}}$,
by a refined set $\mathbb{L} \subset \mathbb{I}_{\text{safe}}$ of perturbations.
We give further discussion and examples of this in subsections \ref{sec:econ} and \ref{sec:biol}.

To proceed we introduce the concept \emph{basin-time}
as the subset of the basin of attraction from which all trajectories returns within a time limit, $t = \tau$.
Based upon the basin-time we build two more resilience measures, $\mathcal{P}^{\tau}$ and $\mathcal{D}^{\tau}$,
by simply replacing basin by basin-time in the constructions of
the nonlocal stability measures $\mathcal{P}$ and $\mathcal{D}$.
To calculate these measures one has to,
for a given $\tau$,
save all initial conditions that returned in time $\tau$ as ``safe" in a set $\mathbb{I}_{\text{safe}}^{\tau}$ and
those that not returned in time $\tau$ as ``unsafe" in a set $\mathbb{I}_{\text{unsafe}}^{\tau}$.
Let $\mathbb{N}_{\text{safe}}^{\tau}$, $\mathbb{N}_{\text{unsafe}}^{\tau}$ and $\mathbb{N}_{\text{tot}}^{\tau}$ denote the safe,
the unsafe and the total number of tested initial conditions,
respectively.
We estimate the resilience measures $\mathcal{P}^{\tau}$ and $\mathcal{D}^{\tau}$ through
\begin{align*}
\mathcal{\hat P}^{\tau} = \frac{\mathbb{N}_{\text{safe}}^{\tau}}{\mathbb{N}_{\text{tot}}^{\tau}} \quad \text{and} \quad
\mathcal{\hat D}^{\tau} = \min_{\substack{x \,\in \,\mathbb{I}_{\text{unsafe}}^{\tau}}} \text{dist}(x, \mathbb{A}).
\end{align*}
In analogue with the stability measures $\mathcal{P}$ and $\mathcal{D}$ we note that
the measure $\mathcal{P}^{\tau}$ is based upon the size of the basin-time while $\mathcal{D}^{\tau}$ is based upon the shape of the basin-time,
and that properties and extensions of the stability measures trivially generalize to these resilience-type measures.
In contrary to the stability measures, 
the resilience measures $\mathcal{P}^{\tau}$ and $\mathcal{D}^{\tau}$ will be applicable also in the setting of simple linear systems having only one global attractor
(in which always $\mathcal{P} = 1$ and $\mathcal{D} =$ constant)
as they record also if the dynamics becomes faster/slower by accounting for the returntime.
The measure $\mathcal{P}^{\tau}$ reflects the probability that the system returns from a random perturbation in time $\tau$,
and $\mathcal{D}^{\tau}$ yields the smallest perturbation from which the system will not return in time $\tau$.




The presented approach involve some degrees of freedom in terms of parameters,
e.g.~the choice of a suitable set of perturbations,
the choice of a distance function for $\mathcal{D}$ and
the choice of time limit $\tau$ for the measures building on the basin-time.
Even though these choices may alter the results,
they should not cause much of a problem since one is always interested in
how a measure reacts on a change of a system,
and thereby a relative change of the value of the measure for fixed parameters.


We will frequently compare our results to the standard local resilience measure given by 
the eigenvalue of the Jacobian matrix at an equilibrium having largest real part. 
We denote this measure by $-\lambda_{\text{max}}$.


\section{Results}
\label{sec:results}

In this section we compare and discuss properties of the suggested nonlocal stability and resilience measures on three different models;
the Solow-Swan model of economic growth, 
a simple electro-mechanical model and 
a stage-structured population model.
Trough all examples we also discuss relations to a local approach for measuring stability and resilience. 
Each model is considered in an independent subsection,
and the reader may start reading the application closest to her/his interest.

\subsection{The Solow-Swan model of economic growth}
\label{sec:econ}

\begin{figure}[h]
\begin{center}
\includegraphics[height = 10cm, width = 13cm]{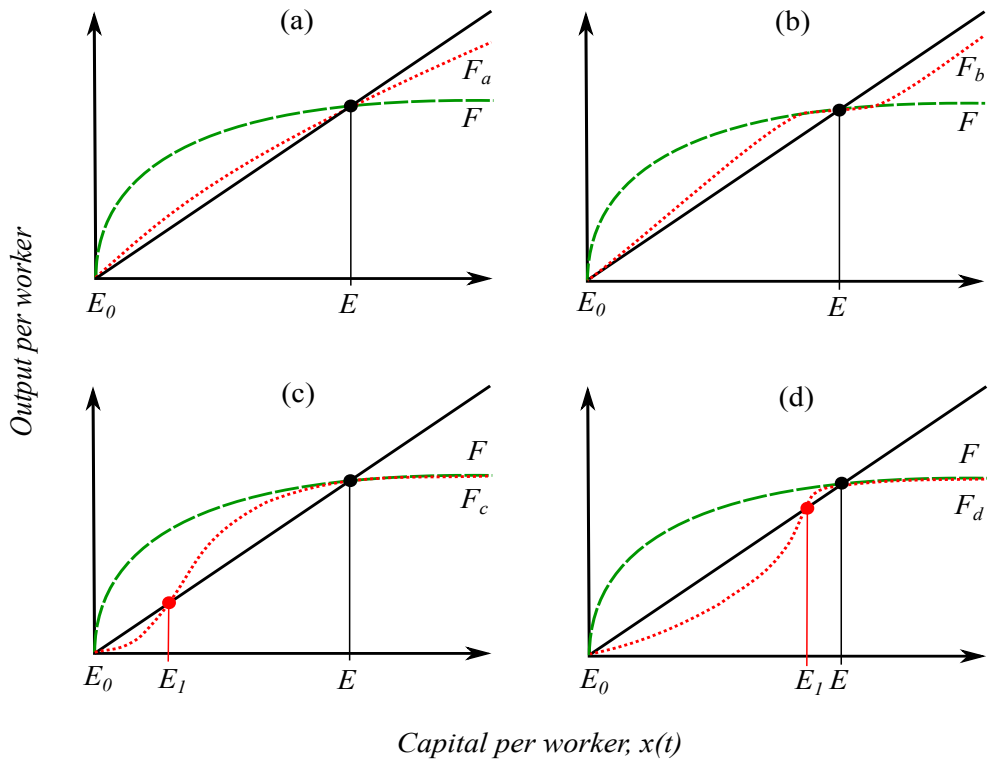} \hspace{1cm}
\caption{
The saved output per worker is given by the function $F$ (green, dashed) of which we consider four variants $F_b$, $F_c$ and $F_d$ (red, dotted).
The maintenance $Cx$ is linear and illustrated by the black solid line.
In (a) and (b) the economy converge to the globally stable equilibrium $E$ given any perturbation.
In (c) and (d) trajectories starting to the left of the unstable equilibrium $E_1$ will reach the origin $E_0$,
while trajectories starting to the right of $E_1$ will reach $E$.}
\label{fig:Solow_model}
\end{center}
\end{figure}

In this section we consider the simple one-dimensional differential equation
\begin{align}\label{eq:solow-swan}
\frac{d x}{dt} \,=\, F(x) \,-\, C x,
\end{align}
where $x(t) \geq 0$, $t > 0$, $C$ is a positive constant and $F(x)$ is a smooth function.
If $x = x(t)$ is capital per worker,
$F(x) = s f(x)$, where $f(x)$ is the production function and $s$ is a constant fraction of the income saved per worker,
and $Cx$ is the investment required to maintain capital per worker,
then the differential equation \eqref{eq:solow-swan} is the well known Solow-Swan model of economic growth (Solow 1956; 
Swan 1956).
With this model, Solow and Swan showed that, under suitable assumptions on $F(x)$, the capital per worker function $x(t)$ reaches an equilibrium after sufficient time.
%
Indeed, the function $F(x)$ is often assumed to grow faster than maintenance $Cx$ as $x$ is small,
and slower than $Cx$ as $x$ is large, see Figure \ref{fig:Solow_model}.
Under such assumption, the Solow-Swan model \eqref{eq:solow-swan} has a globally stable equilibrium ($E$) when saving equals maintenance,
and an unstable equilibrium at the origin:
To the left of $E$, $F(x) > C x$ and therefore $x$ is growing until it reaches $E$,
and to the right of $E$, $F(x) < C x$ and therefore $x$ is decreasing until it reaches $E$.

To test the measures of stability and resilience presented in Section \ref{sec:methods}
we will now consider different variations in the function $F = F(x)$ that make the ``safe" equilibrium $E$ less stable.
Indeed, we consider four variants of perturbing the function $F$,
denoted by $F_a$, $F_b$, $F_c$ and $F_d$,
see Figures \ref{fig:Solow_model} (a)-(d),
and explain how each measure captures
the corresponding consequences of the stability and the resilience in each case.
%

To describe the four cases of stressing the dynamics, we first note that
if $F$ changes to $F_a$ or $F_b$,
as shown in Figures \ref{fig:Solow_model} (a) and (b),
then the Solow-Swan model still has the same equilibria and the global attractor is still equilibrium $E$.
The basin of attraction for $E$ is unchanged in both cases, but
the returntime to $E$ will be longer
(for all perturbations in case $F_a$ and large perturbations in case $F_b$)
since the functions $F_a$ and $F_b$ are closer to the line $Cx$ than $F$ is.
In case of $F_a$ the local curvature near $E$ has changed,
but the curvatures of $F$, $F_b$, $F_c$ and $F_d$ are identical in a small neighborhood of $E$.

In the cases of $F_c$ and $F_d$ the Solow-Swan model has,
besides equilibrium $E$,
a locally stable equilibrium at the origin ($E_0$) and an unstable equilibrium ($E_1$),
see Figures \ref{fig:Solow_model} (c) and (d).
Observe that now $E$ is only \emph{locally} stable, and
the economy can therefore converge to two attractors,
$E_0$ and $E$.
The new equilibrium $E_1$ produces a border of the basins for $E_0$ and $E$;
all trajectories starting to the left of $E_1$ will reach the origin,
while trajectories starting to the right of $E_1$ will reach $E$.

We will focus our discussion on all measures considered in Section \ref{sec:methods},
as well as their relation to local measures.
We do not specify the distance function for $\mathcal{D}$, but the reader may think of it as the Euclidean distance for simplicity.
The exact choice of distribution for perturbations is not important for the below results,
but we assume that perturbations are chosen so that tested initial conditions are distributed with positive probability over the whole interval
displayed in Figure \ref{fig:Solow_model}.
Smaller intervals centered at $E$ can also be considered without altering the main results in this section.

\subsubsection*{Local measures fail to detect changes to $F_b$, $F_c$ and $F_d$}

Since the curvatures of $F$, $F_b$, $F_c$ and $F_d$ coincide in a small neighborhood of $E$,
any local measure, such as e.g.~$-\lambda_{\text{max}}$,
based on the eigenvalues of the Jacobian matrix,
will fail to detect these variations of $F$.
Indeed, the same is true for any function having similar curvature as $F$ near $E$.
This underlines the non-sufficiency of considering only local stability measures;
they will miss all changes of the system except those arbitrary close to $E$.
Moreover, a local analysis of stability and resilience will not reveal the crucial fact that
the modelled economy can converge to two completely different stable states in cases of $F_c$ and $F_d$;
the origin, corresponding to no capital per worker, and $E$.
If $F$ changes to $F_a$, then the local curvature also changes and a local measure will successfully record a decrease in resilience.

\subsubsection*{The measures $\mathcal{P}$ and $\mathcal{D}$ detect changes to $F_c$ and $F_d$ but fail in case of $F_a$ and $F_b$}

The size and the shape of the basin of $E$ will be the same in the three cases of $F$, $F_a$ and $F_b$,
and therefore the measures $\mathcal{P}$ and $\mathcal{D}$ do not work.
Indeed, $\mathcal{D}$ is the distance from $E$ to the origin and $\mathcal{P} = 1$ in all three cases,
even though trajectories will return faster to $E$ in the case of $F$ than with $F_a$ and $F_b$.

If $F$ changes towards $F_d$ through $F_c$,
meaning that the unstable equilibrium $E_1$ moves from the origin towards $E$,
then the size and shape of the basin of $E$ changes and therefore
both measures $\mathcal{P}$ and $\mathcal{D}$
($\mathcal{D}$ is simply the distance between $E$ and $E_1$) will record a decrease in stability,
recall Figure \ref{fig:Solow_model} (c) and (d).
%
In the case of $F_d$, when equilibrium $E_1$ is ``very" close to $E$,
then $\mathcal{D}$ gives the best warning signal that the economy may jump to a completely different stable state (the origin) due to a small perturbation.
The size of basin measure $\mathcal{P}$ will not give such a clear warning signal
since the basin is still large to the right of equilibrium $E$.
Indeed, if we consider in addition that $E_1 \to E$,
then $\mathcal{D} \to 0$ and $\mathcal{P} \to c$ for some constant $c > 0$.

\subsubsection*{The resilience measures $\mathcal{R}$, $\mathcal{R}_{\text{worst}}$, $\mathcal{P}^{\tau}$ and $\mathcal{D}^{\tau}$ detect all considered variations of $F$}

In contrary to a local resilience measure that can only work in case $F_a$,
all suggested nonlocal resilience measures record the crucial changes of the dynamics for all considered variations of $F$.
When $F$ changes to $F_a$ and $F_b$,
the returntime to $E$ will be longer since the functions $F_a$ and $F_b$ are closer to the line $Cx$ than $F$ is.
Therefore, the measures $\mathcal{R}$ and $\mathcal{R}_{\text{worst}}$,
built directly on the returntime, records the destabilization of $E$.
Moreover, in contrary to the basin of equilibrium $E$,
the basin-time of $E$ will,
for any suitable choice of $\tau$,
decrease in size as $F$ changes to $F_a$ or $F_b$.
The reason is that fewer initial conditions will reach $E$ in time $\tau$ when the returntime becomes longer.
Hence, the resilience measures $\mathcal{P}^{\tau}$ and $\mathcal{D}^{\tau}$
also record the destabilization of $E$ as $F$ changes to $F_a$ or $F_b$.

When $F$ changes towards $F_d$ through $F_c$, then,
in addition to the slower convergence of the trajectories towards $E$,
the size of the basin of $E$ also decreases as the new equilibrium $E_1$ comes closer to $E$.
Because $\mathcal{D}^{\tau}$ is bounded from above by the distance from $E_1$ to $E$,
and since $\mathcal{R}_{\text{worst}} \to 0$ immediately if $E_1$ comes into the range of pre-specified perturbations,
these two measures give the clearest warning signals as $E_1$ moves towards $E$.
Measures $\mathcal{R}$ and $\mathcal{P}^{\tau}$ also records a decrease in resilience,
but as they account for that $E$ attracts all trajectories from the right
in a similar way for all three cases $F$, $F_c$ and $F_d$,
they will not approach zero as $E_1 \to E$.
For $\mathcal{R}_{\text{worst}}$ to work properly in this situation
the set of predefined perturbations, $\mathbb{L}$,
should not include perturbations in the whole range from $E$ to the origin,
as that would force the measure down to zero too early as soon as $E_1$ is born.

We have seen that all four nonlocal resilience measures
$\mathcal{R}$, $\mathcal{R}_{\text{worst}}$ , $\mathcal{P}^{\tau}$ and $\mathcal{D}^{\tau}$
are suitable for analysing the economy described by the Solow-Swan model.
Which of these four measures to use depends upon the question to be answered.
The questions related to these measures,
recall Section \ref{sec:intro} and \ref{sec:methods},
should all be relevant for applications in economy.



\subsection{A simple electro-mechanical system}
\label{sec:mech}

\begin{figure}[h]
\begin{center}
\includegraphics[height = 3.5cm, width = 10cm]{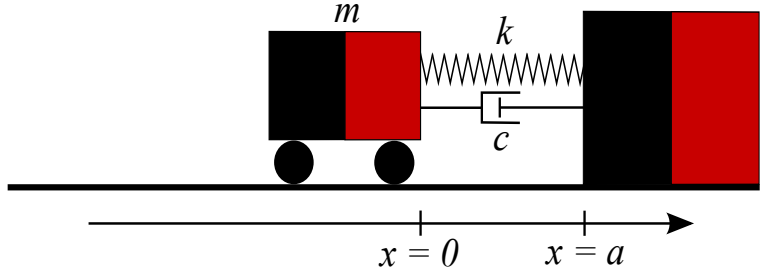}\hspace{2cm}
\caption{
The electro-mechanical machine consists of a wagon with mass $m$ connected to a linear spring with stiffness $k$ and a linear damper with damping factor $c$.
The spring always forces the wagon to $x = 0$ where the spring force is zero.
The wagon is also affected by a nonlinear magnetic force due to a magnet placed at $x = a$.}
\label{fig:wagon}
\end{center}
\end{figure}

In this section we consider an electro-mechanical machine consisting of a wagon with mass $m$,
free to move without friction horizontally in the $x$-direction.
The wagon is connected to a linear spring with stiffness $k$ and a linear damper with damping factor $c$.
The spring always forces the wagon to $x = 0$ where the spring force is zero.
The wagon is also affected by a nonlinear magnetic force from a magnet placed at $x = a$,
see Figure \ref{fig:wagon}.
We assume that the  magnetic force is given by $k_m (x-a)^{-2}$ where the parameter $k_m$ is a magnetic stiffness.
The governing equation of motion is
\begin{align*}
m \,\frac{d^2y}{dt^2} \,=\, - k \, x(t)\, -\, c \frac{dx}{dt} \,+\, \frac{k_m}{(x(t) - a)^2}.
\end{align*}
By introducing another dimension for the velocity,
$y(t) = \frac{dx}{dt}$, we obtain a two-dimensional nonlinear dynamical system on the form \eqref{eq:dyn-syst}:
\begin{align}\label{eq:mek-syst}
\frac{dx}{dt} \,&= \, y(t),\\
\frac{dy}{dt} \,&= \, -\frac{k}{m} \, x(t) \,-\, \frac{c}{m} \, y(t)\, +\, \frac{k_m} {m (x(t) - a)^2},\nonumber
\end{align}
for $x(t) < a$ and $y(t) \in \mathbf{R}^n$.
System \eqref{eq:mek-syst} has (if $k$ is large enough) two equilibria,
of which one is unstable and the other is locally stable ($E$) and constitutes the ``safe" attractor that we will investigate.
There is also always an ``unsafe" attractor at $x = a$ where the wagon slams into the magnet and the machine damages.
More advanced electro-mechanical systems having similar dynamic properties as system \eqref{eq:mek-syst}
arise e.g.~when modelling electric generators and electric motors (see e.g.~Lundstr\"om and Aidanp\"a\"a 2007 and references therein).
One is then typically concerned with at least four-dimensional systems of equations governing the motion of the rotor center
(see e.g. the Jeffcott-rotor model)
and with nonlinear magnetic forces acting between the rotor and the stator.

We will investigate stability and resilience properties of the locally stable equilibrium $E$
as the spring performance decays, using measures suggested in Section \ref{sec:methods}.
In a first case we consider only a decrease of the spring stiffness $k$.
In a second case we assume also that the spring breaks ($k$ is set to zero) if the spring is deformed to quickly,
i.e.~if the wagon moves faster than at a certain speed, $y(t) \geq y_{\text{limit}}$ for some $t > 0$.
In such case, the wagon slams into the magnet at $x = a$ with probability 1.

\subsubsection*{Setting up relevant measures}

In order to specify the set of perturbation to which the machine is to be tested for,
one should consider e.g.~behaviour of connected machines or dynamical systems
as well as other relevant information that can be achieved through the environment of the system.
One may consider only displacements, or only velocity impulses, or combinations of them.
A velocity impulse may be due to a ``smash" on the wagon,
while a displacement may be due to a ``smash" of the fundament of the machine.
The less information you have, the more general set of perturbations may be needed to be considered.
As we have little information in this example, we chose to simply test the machine for perturbations normally distributed in both displacement and velocity
and sample our initial conditions from a two-dimensional normal distribution centered at $E$ with standard deviation 5 in both the $x$- and $y$-direction.
To numerically integrate trajectories from their initial conditions to the attractor,
we used MATLAB's ode-solver ODE45 with standard tolerance settings.
We integrated each trajectory until it either reached the small neighborhood given by a ball of radius 0.01 centered at $E$,
or reached a state where $x(t) > a-0.01$.

As the dynamics of the electro-mechanical system \eqref{eq:mek-syst} is fast and the work of engineers should naturally focus on if the machine damages or not,
rather than the time of recover given a perturbation,
we invoke only our simplest resilience measure $\mathcal{R}$ and put focus on the pure stability measures.
Indeed, as system \eqref{eq:mek-syst} has more than one attractor,
both the size and the shape of basin are crucial to investigate in order to understand stability
and we therefore include
the measures $\mathcal{P}$ and $\mathcal{D}$ in our investigation.
We consider two versions of the latter,
one using simply Euclidean distance ($\mathcal{D}_{\text{Euc}}$) to compare distance to velocities ($x$ to $y$),
and one using an energy-norm to find the distance ($\mathcal{D}_{\text{energy}}$).
Using Euclidean distance may seem naive since it compares velocity to distance in an unmotivated way,
but we chose to include $\mathcal{D}_{\text{Euc}}$ anyway due to its simplicity.

\begin{figure}[h]
\begin{center}
\includegraphics[height = 7.0cm, width = 7.5cm]{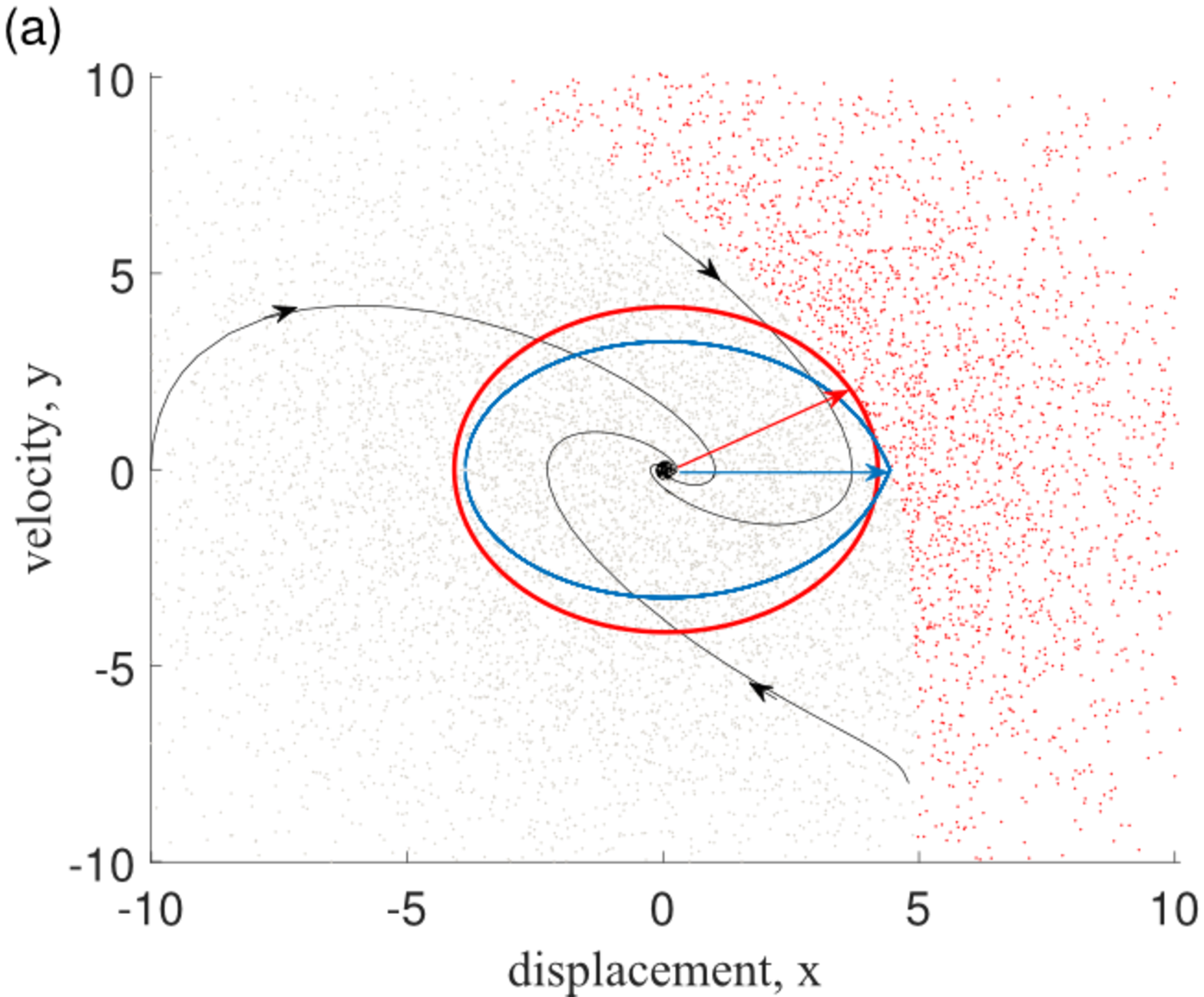}\hspace{1mm}
\includegraphics[height = 7.0cm, width = 7.5cm]{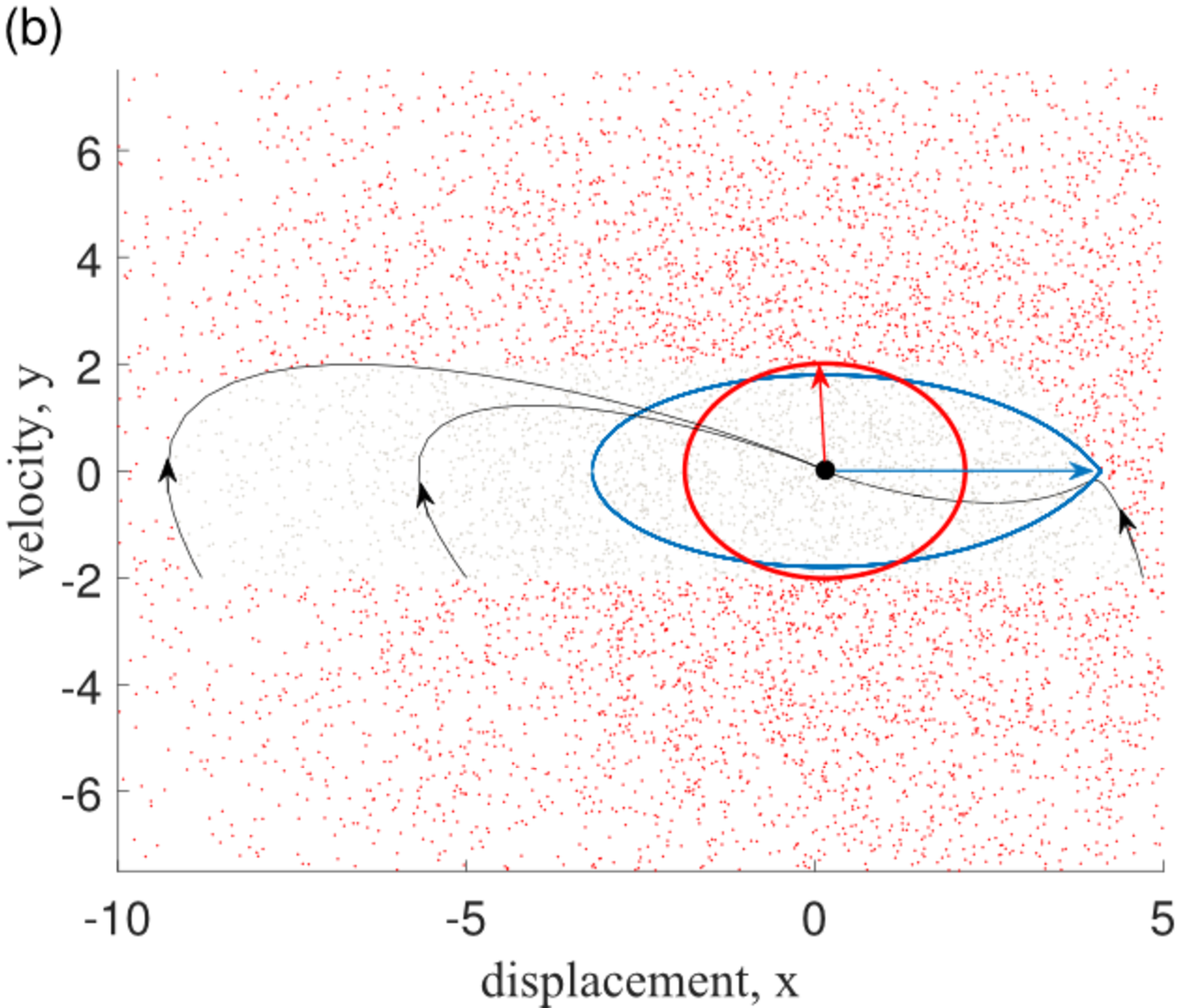}
\caption{
Each dot represents an initial condition drawn from a two-dimensional normal distribution centered at $E$ with standard deviation 5.
A total of 10000 points are drawn in each subfigure.
Grey dots represent safe perturbations as they return to the locally stable equilibrium $E$,
while red dots are unsafe as they do not return to $E$.
The black curves show trajectories of system \eqref{eq:mek-syst} returning to $E$.
The red circle gives an estimate of $\mathcal{D}_{\text{Euc}}$ while the blue level curve of $W(x,y)$ gives an estimate of $\mathcal{D}_{\text{energy}}$.
The red arrows starting at $E$ show the point where the circles tangents the boundary of the basin,
while the blue arrows show the point where the level curves of $W(x,y)$ tangents the boundary of the basin.
%
Parameters are set to $m = c = k_m = 1$ and $a = 5$.
Spring stiffness and speed limit are (a) $k = 0.7$, $y_{\text{limit}} = \infty$
and (b) $k = 0.3$ and $y_{\text{limit}} = 2$.}
\label{fig:basin_mek}
\end{center}
\end{figure}

More natural is the measure $\mathcal{D}_{\text{energy}}$ that will find the least amount of energy needed to push the wagon from the
locally stable equilibrium $E$ into the ``unsafe" attractor at $x=a$.
This is certainly relevant information concerning stability of the machine. 
Giving the wagon an initial velocity $y_0$ implies an input of the kinetic energy of $m\, y_0^2 / 2$.
To find the potential energy input given by an initial displacement $x_0$ we calculate the work needed for moving the wagon
from $E$ to the initial state $x_0$.
We assume that the wagon is moved slowly and that damping is neglected.
We also assume that no energy is transferred back in cases when the wagon moves freely,
due to the magnet force, further away from $E = (x_{eq},0)$.
The potential energy input given by an initial displacement from $x_{eq}$ to $x_0$ is given by
$$
\int_{x_{eq}}^{x_0} \max\left(0, k \, x - \frac{k_m}{(x-a)^2}\right) dx,
$$
and so every initial condition $(x_0, y_0)$ is associated with the energy
\begin{align}\label{eq:energy-norme}
W(x_0,y_0) \, = \, \int_{x_{eq}}^{x_0} \max\left(0, k \, x - \frac{k_m}{(x-a)^2}\right) dx \, + \, \frac{m \,y_0^2}{2}.
\end{align}
We can now define the energy-normed distance and hence our measure $\mathcal{D}_{\text{energy}}$ as
\begin{align}\label{eq:shape_basin_energy}
\mathcal{D}_{\text{energy}} \, = \,  \min_{\substack{x \,\in\, \mathbb{I}_{\text{unsafe}}}} W(x_0,y_0).
\end{align}
While $\mathcal{D}_{\text{Euc}}$ is estimated by the radius of the smallest possible circle, centered at $E$,
not lying entirely in the basin of $E$,
the measure $\mathcal{D}_{\text{energy}}$ is estimated by the smallest value of a constant $C$
such that the level curve $W(x,y) = C$ does not lie entirely in the basin of $E$.
Figure \ref{fig:basin_mek} shows some initial conditions, trajectories as well as the circle
estimating $\mathcal{D}_{\text{Euc}}$ and the level curve estimating $\mathcal{D}_{\text{energy}}$.
Besides giving a stability measure, the construction of $\mathcal{D}_{\text{Euc}}$ and $\mathcal{D}_{\text{energy}}$ reveals
the most critical direction of a perturbation (in the sense of the corresponding measure)
given by the red and blue arrows in Figure \ref{fig:basin_mek}.


\subsubsection*{The measures $\mathcal{D}_{\text{energy}}$ gives the best warning before the machine damages}

\begin{figure}[h]
\begin{center}
\includegraphics[height = 10cm, width = 13cm]{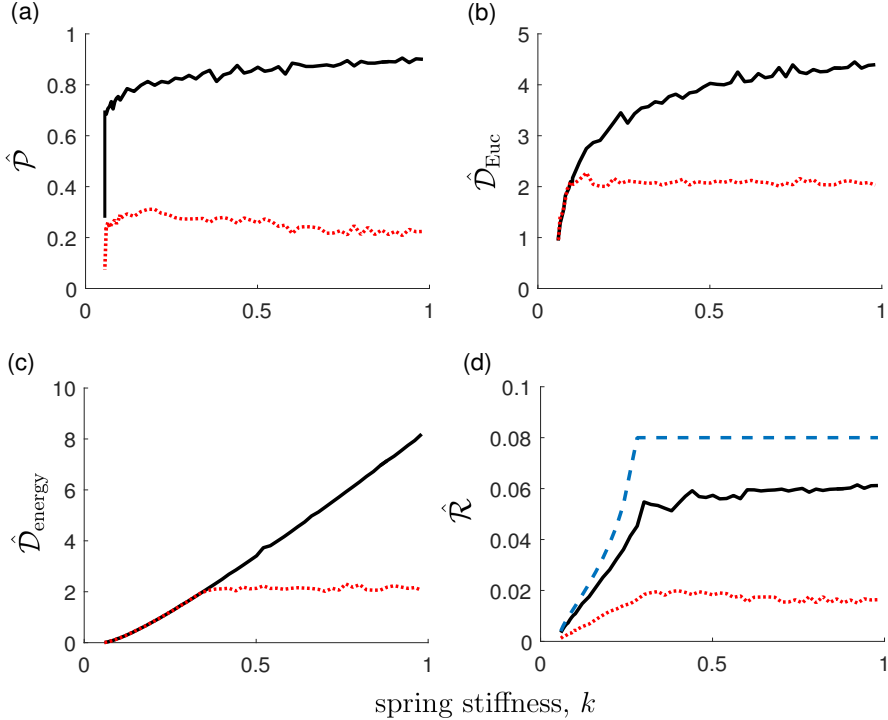}
\caption{
The considered measures as functions of the spring stiffness $k$ in both the case with speed limit (red, dotted) and without speed limit (black, solid).
(a) $\mathcal{P}$,
(b) $\mathcal{D}_{\text{Euc}}$,
(c) $\mathcal{D}_{\text{energy}}$ and
(d) $\mathcal{R}$ and $-\lambda_{\text{max}}\times 0.16$ (blue, dashed).
The energy-based measure of shape of basin, $\mathcal{D}_{\text{energy}}$
gives the clearest warning before the machine damages at $k = 0.054$.
The local measure $-\lambda_{\text{max}}$ can not see the nonlinearity imposed by a speed limit and is therefore the same for both cases.
1000 initial conditions,
randomly sampled from the assumed normal distribution,
are examined for each value of $k$.
Parameters are set to  $m = c = k_m = 1$, $a = 5$ and speed limit $y_{\text{limit}} = 2$.}
\label{fig:measures_mek}
\end{center}
\end{figure}

Figure \ref{fig:measures_mek} shows the measures
$\mathcal{P}$, $\mathcal{D}_{\text{Euc}}$, $\mathcal{D}_{\text{energy}}$, $\mathcal{R}$ and $-\lambda_{\text{max}}$
for decreasing spring stiffness $k$.
At $k \approx 0.054$ the locally stable equilibrium $E$ disappears in a fold bifurcation,
leaving $x = a$ as the only attractor
and the machine damages as the wagon slams into the magnet for any initial condition.
The results for no speed limit is given by the black solid curves
while dotted red curves correspond to the case with speed limit $y_{\text{limit}} = 2$.
Figure \ref{fig:measures_mek} (d) gives also the local measure $-\lambda_{\text{max}}$ (dashed blue),
which, due to its local nature, is unable to record effects of the speed limit and is therefore given by the same curve in both the case with and without speed limit.

The basin of attraction to the safe equilibrium $E$ is relatively large until just before spring stiffness $k$ reaches the critical value of the bifurcation.
Thus, even though $E$ is very close to the boundary of the basin,
so that only a small perturbation to the wagon will damage the machine,
the measure $\mathcal{P}$ will not give a clear warning, see Figure \ref{fig:measures_mek} (a).
This case is similar to Figure \ref{fig:fig_1} (c) in the Introduction;
considering only the size of basin is simply not enough in such situation.
However, one should note that the system is considerably damped,
and by decreasing damping the measure $\mathcal{P}$ would tell more.
This is because less damping allows for more oscillating motions which would force the basin to decrease with decreasing spring stiffness $k$.

As the measure of the shape of basin, $\mathcal{D}_{\text{Euc}}$, records the distance to the boundary of the basin,
it gives a slightly better warning signal, at least in the case without speed limit, c.f. Figure \ref{fig:measures_mek} (a) and (b).
However, in the case with speed limit the measure $\mathcal{D}_{\text{Euc}}$ is, similar to $\mathcal{P}$, nearly constant until just before the crash.
This is because the circle in Figure \ref{fig:basin_mek} (b) is bounded by the speed limit
until just before the spring stiffness $k$ reaches its critical value.

Figure \ref{fig:measures_mek} (c) shows the more natural energy-based measure of the shape of basin, $\mathcal{D}_{\text{energy}}$,
that records the least amount of energy that can push the wagon into the magnet at $x = a$.
This measure approaches zero in a linear way and delivers an early warning signal, in both the cases with and without speed limit, before the crash.
One reason for this is that $\mathcal{D}_{\text{energy}}$ puts focus on the ``most dangerous" direction of perturbation towards the magnet at $x = a$,
in contrary to $\mathcal{D}_{\text{Euc}}$ that uses Euclidean distances;
in terms of energy it is easiest to perturb the wagon in the positive $x$-direction, see the level curves in Figure \ref{fig:basin_mek}.
To stress generality of the measure $\mathcal{D}_{\text{energy}}$,
we point out that even though analytical expressions as \eqref{eq:energy-norme}
for the energy level curves may be difficult to find in general,
numerical calculations of the energy of a displacement should be possible also for complicated systems.

The resilience measures $\mathcal{R}$ and $-\lambda_{\text{max}}$ shown in Figure \ref{fig:measures_mek} (d)
both records a warning before spring stiffness $k$ reaches its critical value and the machine damages.
Observe however that as $\mathcal{R}$ drops significantly as the speed limit is imposed,
the local measure $-\lambda_{\text{max}}$ is unable to take into account the nonlinear effect of the speed limit.
This is because the local behaviour of the system is identical for all possible values of $y_{\text{limit}}$.
This highlights again the unsufficiency of measuring resilience only through local measures;
clearly the system becomes less resilient as $y_{\text{limit}} \to 0$, and a suitable resilience measure should record this fact.

We conclude that, as in case of the Solow-Swan model considered in subsection \ref{sec:econ},
our nonlocal approach is able to record the considered destabilization of the system.
Indeed, we believe that the measure $\mathcal{D}_{\text{energy}}$ is a promising candidate for estimating stability of mechanical systems,
and we further discuss it's generalizations to more general situations of ``small attractors" in Section \ref{sec:discussion}.





\subsection{A stage-structured population model}
\label{sec:biol}

In this section we consider a stage-structured consumer-resource population model introduced by
de Roos et al. (2008). 
The model is based upon the assumption that individuals are composed into two stages;
juveniles and adults, depending only on their size.
Juveniles are born with size $s_{\text{born}}$ and grow until
they reach the size $s_{\text{max}}$ at which they cease to grow,
mature, and become adults.
Juveniles use all available energy for growth and maturation,
while adults do not grow and instead invest all their energy in reproduction.
The total juvenile biomass is denoted by $J = J(t)$ while the total adult biomass is denoted by $A = A(t)$.
Both juveniles and adults foraging on a shared resource $R = R(t)$.
The net biomass production per unit biomass for juveniles and adults
equals the balance between ingestion and maintenance rate $T$ according to
\begin{equation*}
w_{\text{J}}(R)\,=\, \max\left\{ 0,\, \sigma I_{\max} \frac{R}{H+R} - T \right\} \quad \text{and} \quad
w_{\text{A}}(R)\,=\, \max\left\{ 0, \, \sigma q I_{\max} \frac{R}{H+R} - T \right\},
\end{equation*}
where $q$ describes the difference in ingestion rates between juveniles and adults,
$\sigma$ represents the efficiency of resource ingestion,
$H$ is half-saturation constant of the consumers,
and the maximum juvenile and adult ingestion rates per unit biomass equal $I_{\max}$ and $qI_{\max}$, respectively.
The model incorporates stage-selective harvesting by allowing for separate harvesting rates
on juveniles $(h_{\text{J}})$ and on adults $(h_{\text{A}})$,
in addition to the background mortality rates $d_{\text{J}}$ and $d_{\text{A}}$.
Resource turn-over rate is denoted by $r$ and $R_{\max}$ is the maximum resource density.
The model consists of the following three-dimensional dynamical system:
\begin{eqnarray}\label{eq:fish-syst}
\frac{dJ}{dt} &=& \left( w_{\text{J}}(R) - v(w_{\text{J}}(R)) - d_{\text{J}} - h_{\text{J}}\right) J \,+\, w_{\text{A}}(R)\,A,\nonumber \\
\frac{dA}{dt} &=& v(w_{\text{J}}(R))J \,-\, \left( d_{\text{A}} + h_{\text{A}}\right) A,\\
\frac{dR}{dt} &=& r(R_{\max}-R)\nonumber \,-\, I_{\max}\frac{R}{H+R}\left(J+qA\right),\nonumber
\end{eqnarray}
where $J(t),\, A(t),\, R(t) \in \mathbf{R}_+$ and
\begin{align*}
v(x) \,=\, \frac{x - d_{\text{J}} - h_{\text{J}}} {1-\left(s_{\text{born}}/s_{\text{max}}\right)^{1-(d_{\text{J}} + h_{\text{J}})/x}},
\end{align*}
for $x \neq d_{\text{J}} + h_{\text{J}}$
and $v(d_{\text{J}}+h_{\text{J}}) = -(d_{\text{J}} + h_{\text{J}})/\log(s_{\text{born}}/s_{\text{max}})$.
The function $v(x)$ describes the maturation rate by determining how fast juvenile biomass is transferred into adult biomass.
We adopt the parameter values $H = T = r = 1$, $I_{\text{max}} = 10$, $d_{\text{J}} = d_{\text{A}} = 0.1$,
$q = 0.85$,
$\sigma = 0.5$,
$R_{\max}= 2$,
$s_{\text{born}} / s_{\text{max}} = 0.01$
(Meng et al. 2013; de Roos et al. 2008). 

Concerning the dynamics of the model,
extensive numerical investigations indicate that there exists always a unique global attractor,
and,
depending on the parameter values,
this attractor is either a positive equilibrium
$(J, A, R) = (J_{\text{eq}}, A_{\text{eq}}, R_{\text{eq}})$
or an extinction equilibrium given by $(J, A, R) = (0, 0, R_{\text{max}})$
(Meng et al. 2013). 
This means that the basin of attraction is always the whole $\mathbf{R}^3_+$
and therefore the size and shape of basin will not deliver any information.
In fact, we always have $\mathcal{P} = 1$ and $\mathcal{D} = $ constant.
However, as we will see in the following subsections,
we will obtain the sought after information by recording the returntime through
the resilience measures given in Section \ref{sec:methods}.

Our first aim here is to show how to use the resilience measures presented in Section \ref{sec:methods}
to investigate the consequences of the populations ability to recover the globally stable equilibrium,
given a perturbation,
as harvest pressure increases.
Our second aim for this section is to show how to apply the same measures of resilience in order to compare
the efficiency of different harvesting strategies when accounting for both yield and resilience.

\subsubsection*{Setting up relevant measures}

For this application we consider perturbations of sizes relative to the actual
size of the population biomass at equilibrium.
This choice is mainly based on our believes that it is reasonably that
a fraction, rather than an absolute amount,
of the population is eliminated through e.g.~illegal harvesting,
tough weather conditions, sudden diseases or due to some other reason.
Further motivations for considering perturbations of relative size is given by the fact that the populations equilibrium biomass can vary drastically as function
of the harvesting rates and other parameters,
making absolute measures nonrelevant.

We will consider two approaches of measuring the resilience of
system \eqref{eq:fish-syst} when harvesting pressure increases.
In a first case \emph{including resource} we test the system for perturbations
in all three dimensions; i.e.~juveniles, adults and the resource.
This approach is relevant when studying how the system reacts when both the population and the resource
are perturbed, e.g.~by sudden environmental changes.
In a second case \emph{excluding resource} we leave the resource at equilibrium and
test the system for perturbations only in the two dimensions of juveniles and adults.
This is relevant when studying how the system reacts when the population,
but not (explicitly) the resource,
is perturbed e.g.~by sudden diseases, by illegal harvesting or
by errors in an implementation of a harvesting strategy.
We proceed by specifying the four resilience measures
$\mathcal{R}$, $\mathcal{R}_{\text{worst}}$, $\mathcal{P}^{\tau}$ and $\mathcal{D}^{\tau}$,
given in Section \ref{sec:methods},
for the approach including resource as well as for the approach excluding resource.

\textit{Including resource.}
We impose perturbations randomly by drawing initial conditions from a three-dimensional normal distribution
centered at the equilibrium $(J_{\text{eq}}, A_{\text{eq}}, R_{\text{eq}})$
with standard deviation given by
$\frac12\left( J_{\text{eq}}, A_{\text{eq}}, R_{\text{eq}} \right)$ in each dimension.
Initial conditions including non-positive values are excluded.

For the measures $\mathcal{P}^{\tau}$ and $\mathcal{D}^{\tau}$ based on the size and the shape of the basin-time
we have to specify $\tau$.
In general, to find a suitable value of $\tau$ one should recall the sought after information and that
$\mathcal{P}^{\tau}$ estimates the probability that the population returns to equilibrium
before time $\tau$ given a random perturbation,
while $\mathcal{D}^{\tau}$ estimates the largest perturbation
that can be given to the population if it is demanded to return before time $\tau$.
The value of $\tau$ should also be in the range of typical returntimes of the system,
and for our application we will use $\tau = 5$ for both measures.
To find a relevant version of the shape of basin-time measure, $\mathcal{D}^{\tau}$,
we follow the above reasoning of relative distances to the population size and define
\begin{align}\label{eq:D-tau}
\mathcal{D}_{\text{rel}}^{\tau} \,\approx\,  \min_{\substack{x \,\in\, \mathbb{I}_{\text{unsafe}}^{\tau}}}
\sqrt{ \left( \frac{x_1 - J_{\text{eq}}}{J_{\text{eq}}}\right)^2 + \, \left(\frac{x_2 - A_{\text{eq}}}{A_{\text{eq}}}\right)^2 + \, \left(\frac{x_3 - R_{\text{eq}}}{R_{\text{eq}}}\right)^2}.
\end{align}

For the resilience measure $\mathcal{R}_{\text{worst}}$,
estimating ``worst resilience using the slowest recovery,"
we use a smaller set of perturbations ($\mathbb{L}$) than
for the other measures.
In particular, we restrict perturbations to a smaller neighborhood of the attractor
by taking $\mathbb{L}$ as the subset of the above defined initial conditions $x = (x_1, x_2, x_3)$ satisfying
\begin{align}\label{eq:ellipse}
\left(\frac{x_1 - J_{\text{eq}}}{J_{\text{eq}}}\right)^2 + \, \left(\frac{x_2 - A_{\text{eq}}}{A_{\text{eq}}}\right)^2
+ \, \left(\frac{x_3 - R_{\text{eq}}}{R_{\text{eq}}}\right)^2 \, = \, \frac{1}{4}.
\end{align}
Doing so we ensure that a perturbation for $\mathcal{R}_{\text{worst}}$ never removes more than half the juvenile biomass,
or half the adult biomass,
from the population;
indeed, if $J_{\text{eq}} = A_{\text{eq}} = R_{\text{eq}} = 1$
then \eqref{eq:ellipse} is a circle centered at $(1,1,1)$ with radius $\frac{1}{2}$.

\textit{Excluding resource.}
In this case initial conditions are drawn from a two-dimensional normal distribution
centered at $(J_{\text{eq}}, A_{\text{eq}})$
with standard deviation $\frac12\left( J_{\text{eq}}, A_{\text{eq}}\right)$.
The measure $\mathcal{D}^{\tau}$ is defined through \eqref{eq:D-tau} but without the third term under the square rot.
Similarly, the set $\mathbb{L}$ for $\mathcal{R}_{\text{worst}}$ is defined as those initial conditions that end up in the ellipse given by
\eqref{eq:ellipse} but without the resource term.

\begin{figure}[h]
\begin{center}
\includegraphics[height = 6.5cm, width = 7.5cm]{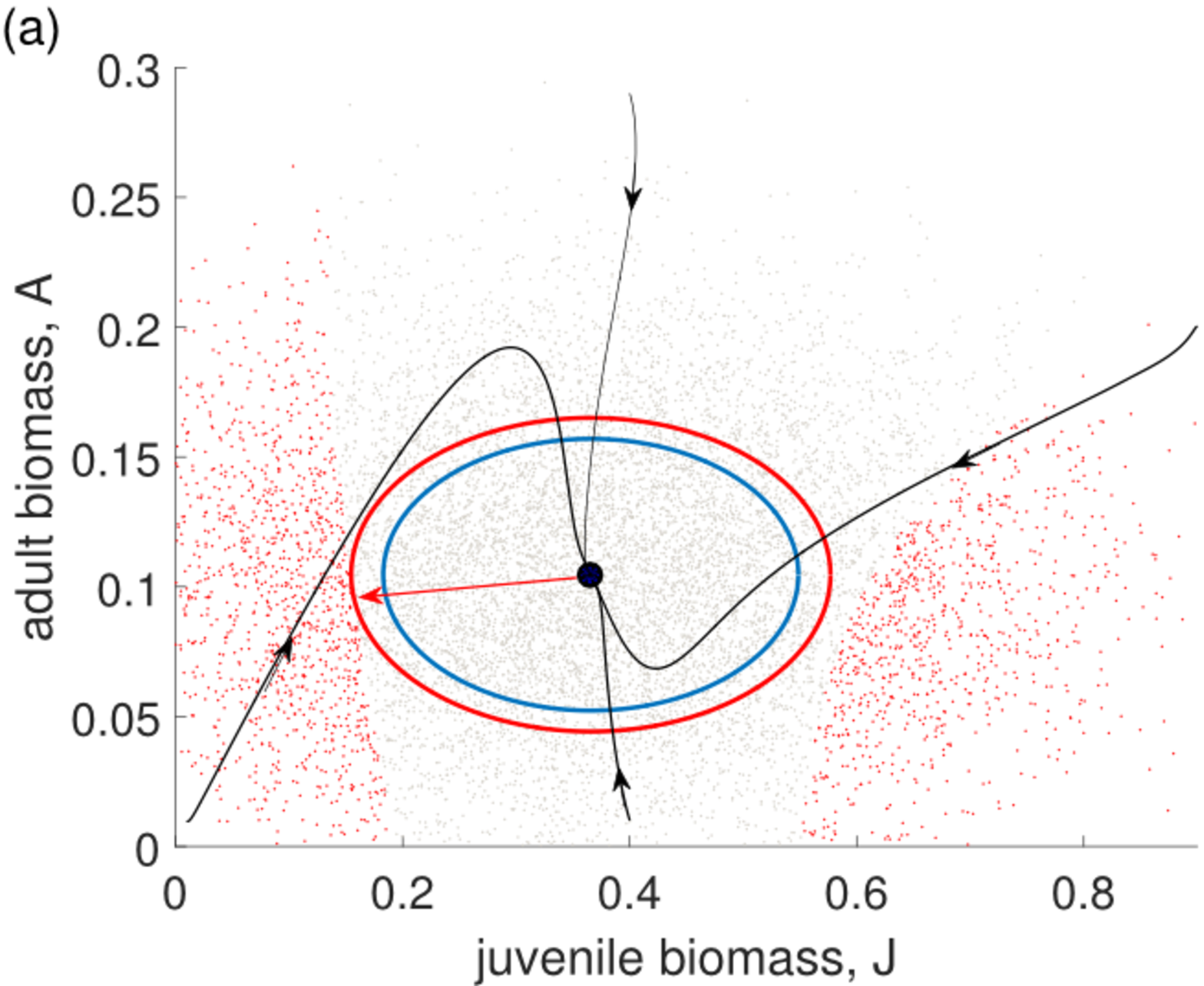}
\includegraphics[height = 6.5cm, width = 7.5cm]{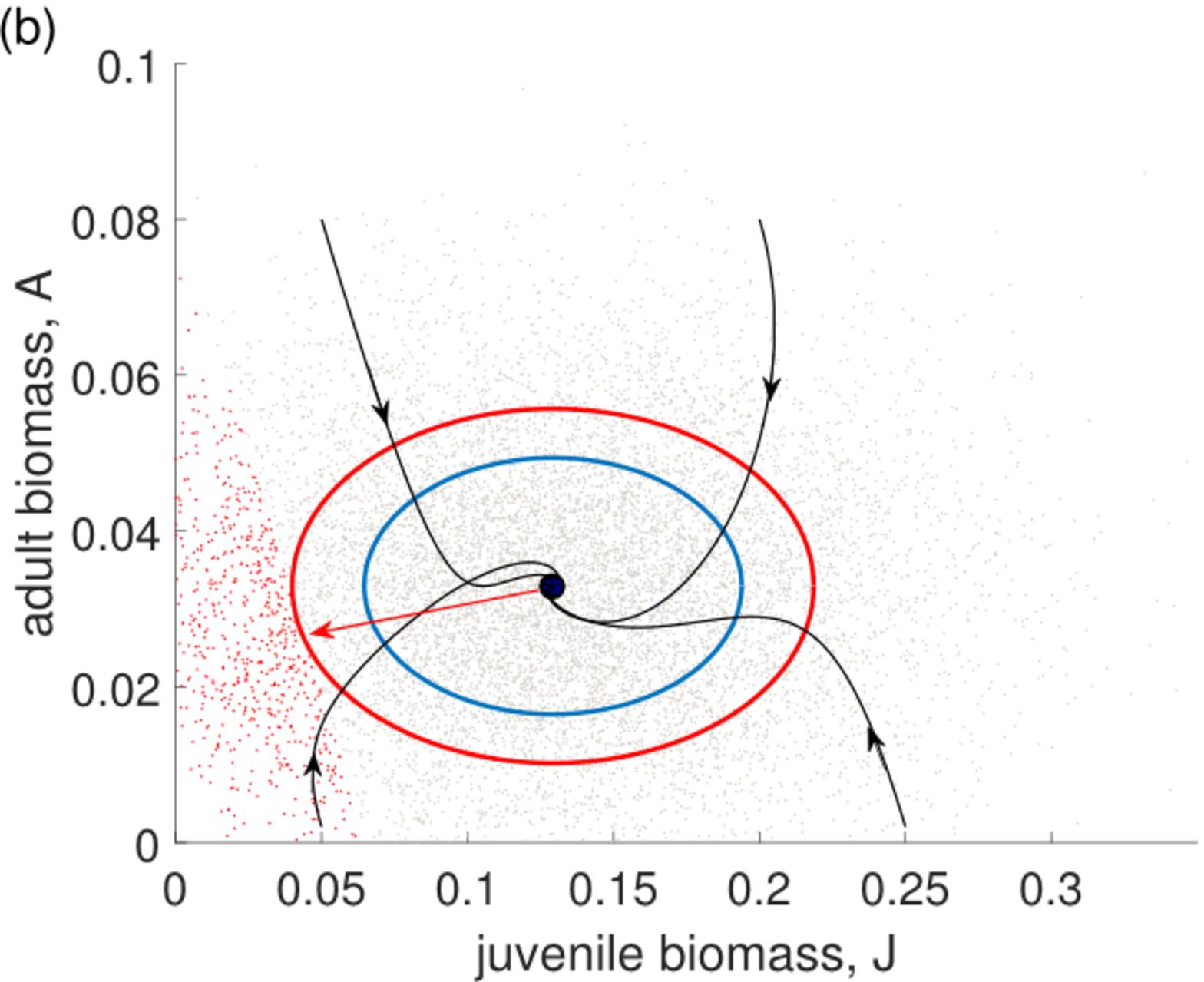}
\caption{
Each dot represents a tested initial condition drawn from the two-dimensional normal distribution centered at the equilibrium (large black dot).
Grey dots forms the set basin-time as corresponding trajectories return to the globally stable equilibrium before time $\tau$,
while trajectories starting at red dots have a returntimes that exceed $\tau$.
The black curves show trajectories of system \eqref{eq:fish-syst} returning to the equilibrium.
The red ellipses give estimates of the resilience measure $\mathcal{D}_{\text{rel}}^{\tau}$ while the blue ellipses define perturbations for the worst-case measure $\mathcal{R}_{\text{worst}}$.
The red arrows give the direction of the shortest distance (in the sense of $\mathcal{D}_{\text{rel}}^{\tau}$) to the boundary of the set basin-time.
A total of 10000 perturbations are tested in each subfigure, $\tau = 5$ and
harvest pressure is (a) $h_j = h_a = 0.5$ and (b) $h_j = h_a = 1.5$.
}
\label{fig:fish-basin}
\end{center}
\end{figure}

\noindent
Figure \ref{fig:fish-basin} shows the basin-time,
some trajectories,
the ellipse \eqref{eq:ellipse} defining the set $\mathbb{L}$ (blue),
as well as the ellipse giving the estimate for the size of basin-time measure $\mathcal{D}_{\text{rel}}^{\tau}$ (red)
in the case excluding resource.

To numerically integrate trajectories from their initial conditions to the attractor,
we used MATLAB's ode-solver ODE45 with standard tolerance settings.
In both the case including resource and excluding resource,
we integrated each trajectory until it either reached the small neighborhood given by the ellipse \eqref{eq:ellipse} but with right hand side 0.01.

\begin{figure}[h]
\begin{center}
\includegraphics[height = 10cm, width = 13cm]{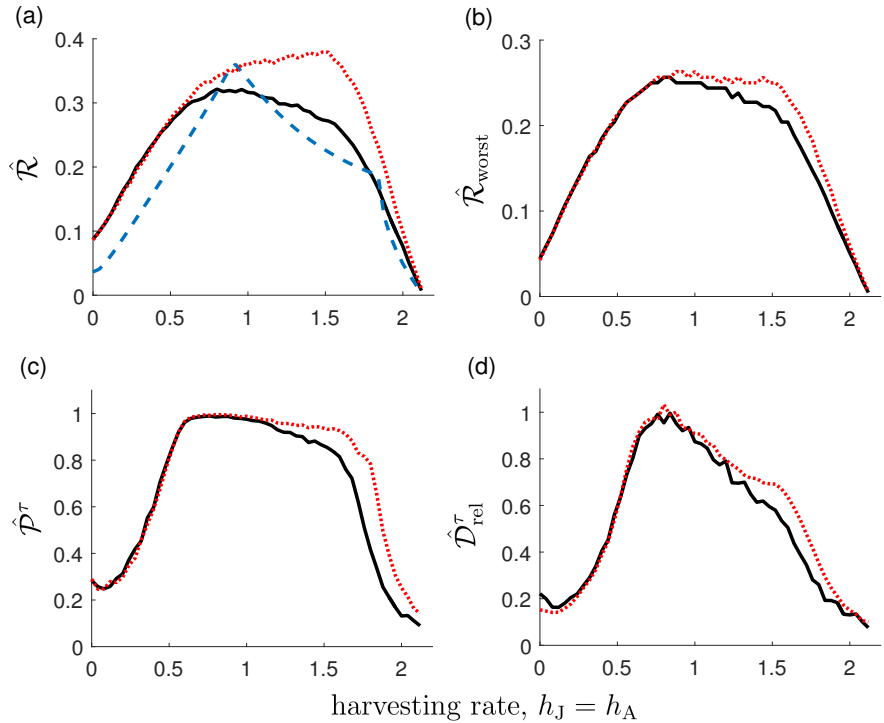}
\caption{
The considered measures as functions of harvest rates $h_{\text{J}} = h_{\text{A}}$ in both the case with resource (black, solid) and without resource (red, dotted).
(a) $\mathcal{R}$ and $-\lambda_{\text{max}} \times \frac{1}{3}$ (blue, dashed),
(b) $\mathcal{R}_{\text{worst}}$,
(c) $\mathcal{P}^{\tau}$ and
(d) $\mathcal{D}_{\text{rel}}^{\tau}$.
All five measures agree on similar results.
2000 initial conditions are examined for each value of $h_{\text{J}} = h_{\text{A}}$.
}
\label{fig:measures_fish}
\end{center}
\end{figure}

\subsubsection*{All resilience measures $\mathcal{R}$, $\mathcal{P}^{\tau}$, $\mathcal{D}_{\text{rel}}^{\tau}$ and $\mathcal{R}_{\text{worst}}$ give similar results}

Figure \ref{fig:measures_fish} shows the suggested measures
$\mathcal{P}^{\tau}$, $\mathcal{D}_{\text{rel}}^{\tau}$, $\mathcal{R}$ and $\mathcal{R}_{\text{worst}}$
together with $-\lambda_{\text{max}}$ for increasing harvesting pressure.
We assume here equal harvesting rates on juveniles and adults, i.e. $h_{\text{J}} = h_{\text{A}}$.
The case with resource is given by black solid curves,
while dotted red curves correspond to the case without resource.
In Figure \ref{fig:measures_fish} (a),
the dashed blue curve gives the local resilience measure $-\lambda_{\text{max}}$.

All five measures agree on a similar result;
the resilience first increases with harvesting pressure,
up to a maximum at approximately $h_{\text{J}} = h_{\text{A}} = 0.8$,
after which it decreases to its lowest level before the population goes extinct.
All five measures of resilience behave as
similar ``smooth" functions of the harvest rates,
which indicates that system \eqref{eq:fish-syst} has a ``weak" nonlinearity;
the linearized system probably agrees well with the original one in a relatively large neighborhood of the investigated equilibrium.

When comparing the approaches of including and excluding the resource in the measures constructions,
we see that the curves corresponding to the approach including resource are lower than the curves for resource excluded.
This is true mainly for relatively high harvest rates and for the measure  $\mathcal{R}$,
see Figure \ref{fig:measures_fish} (a).
Thus, excluding resource will in this situation be in favour of a slightly higher harvest pressure.

\begin{figure}[h]
\begin{center}
\includegraphics[height = 11cm, width = 13cm]{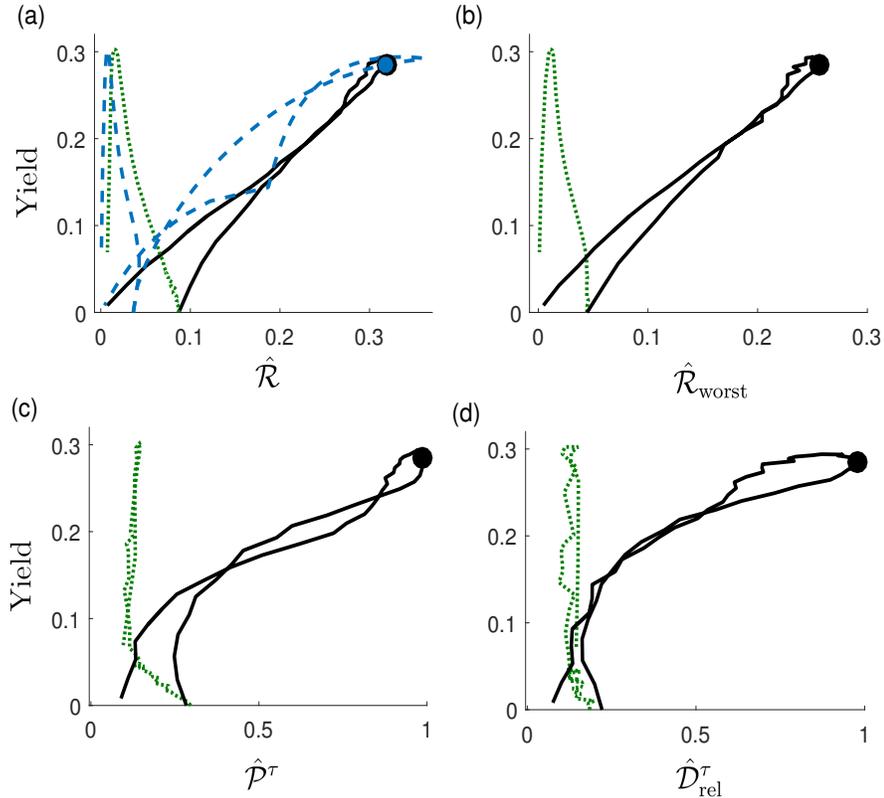}
\caption{
The yield as functions of the five measures of resilience for increasing harvesting pressure. Equal harvesting (black, solid) and adult harvesting (green, dotted).
(a) $\mathcal{R}$ and $-\lambda_{\text{max}} \times \frac{1}{3}$ (blue, dashed),
(b) $\mathcal{R}_{\text{worst}}$,
(c) $\mathcal{P}^{\tau}$ and
(d) $\mathcal{D}_{\text{rel}}^{\tau}$.
The dots corresponds to the strategy $h_{\text{J}} = h_{\text{A}} = 0.8$ which seems to be one of the best possible strategies.
In this simulation, 2000 initial conditions are examined for each harvesting rate combination.
We consider here the case including the resource in the measure constructions.}
\label{fig:pareto_fish}
\end{center}
\end{figure}

\subsubsection*{Comparing harvesting strategies using the measures $\mathcal{R}$, $\mathcal{R}_{\text{worst}}$, $\mathcal{P}^{\tau}$ and $\mathcal{D}_{\text{rel}}^{\tau}$}
\label{sec:biol_2}

In this section we 
use the proposed measures of resilience to compare the strategy of harvesting only adult biomass (adult harvesting),
$h_{\text{J}} = 0, h_{\text{A}} > 0$,
to harvesting both juvenile and adult biomass at the same rates (equal harvesting), $h_{\text{J}} = h_{\text{A}} > 0$.
For a given harvesting strategy we can find the yield (catch) as
\begin{align*}
\text{Yield} \left( h_{\text{J}}, h_{\text{A}} \right) \, = \, h_{\text{J}} J_{\text{eq}} \,+\,  h_{\text{A}} A_{\text{eq}}.
\end{align*}
To find out which strategy that is most efficient when accounting for both yield and resilience,
we plot the yield as functions of the five resilience measures
$\mathcal{P}^{\tau}$, $\mathcal{D}_{\text{rel}}^{\tau}$, $\mathcal{R}$, $\mathcal{R}_{\text{worst}}$
and $-\lambda_{\text{max}}$ in Figure \ref{fig:pareto_fish}.
It is clear, in terms of all resilience measures,
that equal harvesting performs better than adults harvesting.
This is true since for any given yield (close to the maximum yield that can be obtained)
equal harvesting gives much higher resilience,
in terms of all measures, than adult harvesting.
This result strengthen previous results by Lundstr\"om et al. (2016) 
arguing for equal harvesting.
In particular, in that paper the authors evaluate harvesting strategies with respect to both yield and conservation using two population models, of which one is model \eqref{eq:fish-syst}.
They use four measures to quantify conservation,
and their resilience measure is similar to the presented measure $\mathcal{R}$.
%

We have seen that all nonlocal resilience measures
$\mathcal{R}$, $\mathcal{R}_{\text{worst}}$ , $\mathcal{P}^{\tau}$ and $\mathcal{D}^{\tau}$,
as well as the local measure $-\lambda_{\text{max}}$,
give similar results in the sense of recommending the same harvesting strategy.
This does not mean that the presented nonlocal approach is superfluous.
Indeed, it tells that the system behaves ``smooth and linearly" near the attractor and that
we may not see sudden jumps or other unexpected behaviour in the dynamics.
This is valuable information that can not be obtained through a local measure based on a linearized system.
Which of the suggested nonlocal measures to include in a study depends upon the question to be answered.
All questions related to these measures (see Section \ref{sec:intro} and \ref{sec:methods})
should be relevant for a variety of applications in biology.



\setcounter{equation}{0} \setcounter{theorem}{0}

\section{Discussion}
\label{sec:discussion}


We have considered two nonlocal stability measures 
and four nonlocal resilience measures based on properties of the systems basins,
on the introduced concept basin-time, as well as on the returntime to the attractor.
We showed how to quantify stability and resilience using these measures on three different dynamical systems modelling economy,
electro-mechanics and stage-structured populations.
We compared and discussed our results in relation to local measures through the paper,
and explained which measures that give what information in each case.
We concluded that, in contrary to a local approach,
our measures give a good understanding of the stability properties of the attractor in all examples considered.
Due to this fact and due to its simplicity,
we believe that the presented methods have a large potential to increase the understanding of stability 
in a wide range of research fields involving applied dynamical systems.
We underline that an implementation of our approach reveals,
with a high probability,
if there is another unexpected attractor in the range of the considered perturbations.
Indeed, as our method implies testing the system for a large number of initial conditions
even hidden attractors (see e.g.~Dudkowski et al. 2016) can automatically be found
(when producing Figures \ref{fig:measures_mek} and \ref{fig:measures_fish} more than 50000 initial conditions where tested in each system).

When deciding which measures to include in a stability analysis we first note that
measuring the size and the shape of the basin should always be done in some way when dealing with nonlinear systems.
Thus, if no analytical estimates of the basin is available,
measures $\mathcal{P}$ and a suitable version of $\mathcal{D}$ should always be relevant.
In case the returntime is important then a resilience measure should be included.
Which of the resilience measures to use depends on the question to be answered;
you should choice the measure that best give the answer to what you want to know about the systems behaviour.
We recall that calculating all suggested measures is not much more expensive than calculating one measure,
and the more measure you include, the more information of the dynamics you get.

It seems difficult to construct an overall ``best measure" in the general setting of a wide range of applications,
but it is of course possible to build integrated measures that records the relevant information in more general settings than our suggested measures.
However, that may ruin the simplicity in the measures construction and thereby the results will also be harder to interpret.
Therefore, we choice to stay with simple measure constructions
and recommend instead to use several simple measures simultaneously.




\subsubsection*{``Small" attractors can be handled simply as an equilibrium}

We have demonstrated our approach on three models in which the stability and the resilience were investigated on the simplest attractor,
an equilibrium.
The presented ideas are, however,
applicable also in cases of general attractors, not only in the setting of an equilibrium.
In such general case it may be more difficult to implement the ideas as it is usually not easy
e.g.~to calculate the distance from the attractor to the boundary of the basin,
and to determine when a trajectory has returned to the attractor.

However, all the suggested measures are simple to calculate for some complicated attractors.
Indeed, if the attractor is ``small" in the sense that it can be contained in a ball $B(y,r)$,
with radius $r$ and center $y$,
where $r$ is small in comparison to the perturbations one wish to test the system for,
then one can use a similar procedure as if the attractor was an equilibrium at $y$.
This is independent of the type of the attractor; it can be periodic, quasi-periodic or even chaotic.
To explain this,
recall that when the attractor is an equilibrium at $y$,
then we may use a ball $B(y,r)$, $r$ very small, to define the neighborhood
determining when a trajectory has returned to the equilibrium and the simulation can be stopped.
In the setting of an arbitrary attractor contained in a ball $B(y,r)$,
we can proceed similarly;
when the trajectory has entered $B(y,r)$ and stayed there a ``sufficient" time,
then we can assume that the attractor has been reached.
The returntime is thereafter determined through the moment when the trajectory entered the ball $B(y,r)$ for the last time.

This is an important point 
because these situations
occur frequently in applications of dynamical systems.
Indeed, equilibria often bifurcate into ``small" but complicated attractors due to variations in model parameter values.
To explain a typical such situation, consider a ``perfect" electric generator in which the rotor is balanced and no other external forcing on the rotor is present.
In such situation the rotor rotates at a stable equilibrium.
However, in reality there are always imperfections such as e.g.~mass imbalance on the rotor, shape deviations of the rotor and the stator 
as well as external forces transferred from nearby machines such as turbines.
When adding the effects from such phenomena to the perfect rotor model,
the stable equilibrium will most probably bifurcate into more complicated dynamic behaviour,
e.g.~mass unbalance implies simple periodic motion.
Engineers have to keep track of these effects when designing the machine.
Typically, large amplitudes are unwanted and the generator should be constructed so that the attractor is ``small" in the above sense,
i.e.~so that the rotor stays close to the previously existing equilibrium (the attractor in case of the perfect machine) under operation.
As our approach can handle these situations in a simple way,
we believe that the suggested measures will constitute useful tools
for engineers when designing new generators and motors as well as when performing restore work and understanding behaviour of existing machines.
Simple versions of $\mathcal{P}$ and $\mathcal{D}$ have already been used in this manner,
see Lundstr\"om and Aidanp\"a\"a (2007). 


\subsubsection*{Alternative measures and topics for future research}

We have considered two nonlocal stability measures and four nonlocal resilience measures.
As we believe in a value of straightforward interpretations of a measures result in terms of application,
our choice of suggested measures rely heavily on the simplicity in the construction.
Anyway, we choice to define the resilience measures $\mathcal{R}$ and $\mathcal{R}_{\text{worst}}$
through the reciprocal of the returntime, see Section \ref{sec:methods}.
An even more straightforward way to record the systems ability to return fast would be to consider the returntime explicitly.
For bistable systems however, one then has to integrate the effect of the ``unsafe" perturbations, of which trajectories not return to the attractor,
in suitable way.
By taking the reciprocal of the returntime and consider resilience as return-rate,
we solved this by giving zero resilience to unsafe perturbations.
Another reason for considering the resilience, in place of the returntime explicitly,
is the more direct connection to the literature and established resilience measures such as the largest eigenvalue.
In connection to this discussion,
we mention that the function $1 / (T(x) + t_{\epsilon})$,
used in the definitions of $\mathcal{R}$ and $\mathcal{R}_{\text{worst}}$,
may be replaced by any decreasing function such as e.g.~$e^{-T(x)}$ (which decreases faster to zero and therefore punish more for long returntimes).
If the attractor is globally stable, one may also take the mean of returntimes before taking the reciprocal (Lundstr\"om et al. 2016).

There are several suggestions on alternatives to standard resilience measures in the literature.
Neubert and Caswell (1997) consider several alternative measures (e.g.~reactivity and maximum possible amplification)
to resilience together with local methods for calculating them.
They demonstrate their measures on several examples and discuss also disadvantages of local methods for calculating stability and resilience measures.
Extensions of the work of Neubert and Caswell has been done by e.g.~Arnoldi et al. (2016) and Arnoldi et al. (2017).
These papers consider small perturbations and put focus on local measures.
Mitra et al. (2015) 
consider an integrative quantifier of resilience based on both local and nonlocal theory.
Their measure builds upon the three aspects of of resilience by Walker et al. (2004); latitude, resistance and precariousness.
Hellmann et al. (2016) 
consider the nonlocal measure \emph{survivability}, reflecting
the likelihood that the trajectory of a random initial condition leaves a certain desired (or allowed) subset of the systems state-space.
Interesting topics for future research is (1) to consider nonlocal versions of existing local measures and
(2) to compare and discuss relations between existing measures and those considered here.

Even though our approach seems relatively fast
(the simulations producing Figures \ref{fig:measures_mek} and \ref{fig:measures_fish} finished within approximately one hour on a standard laptop computer)
it would be valuable to make it faster.
Indeed, the dynamical systems considered in our examples are simple and low-dimensional.
More advanced models will result in longer computational time.
When using random sampling of initial conditions in order to impose the perturbations on the system,
an efficient variance reduction technique will be useful for reducing the computational time.
Indeed, we have to sample a relatively large number of initial conditions from a probability distribution in order to
get small enough variance of our estimates, and each sampled trajectory is expensive in terms of computational time.
A variance reduction method will ensure that the sampled initial conditions are well spread over the assumed probability distribution.
This reduces the variance of the estimate and therefore less initial conditions need to be tested,--
resulting in less computational time.
We intend to develop a variance reduction method,
suitable for the presented measures,
in an upcoming paper.

Another interesting topic for future research would be to create an efficient approach for finding estimates
of the neighborhood to attractors in nonlinear systems in which the
linearized system approximates the original system well.
One could e.g.~classify dynamical systems based on the nature of their nonlinearity so that ``weak" means
that the linearization holds in a large neighborhood,
and ``strong" means that it holds only in a small neighborhood.
Such results would complement the linear local stability approach by giving knowledge about the magnitude
of perturbations that can be handled by the already well established
local stability and resilience measures.


\subsubsection*{Concluding words}

We have presented two nonlocal stability measures and four nonlocal resilience measures
in order to complement the widespread local stability-measures approach based on linearizations.
The measures are simple and easy to implement on a standard laptop computer.
All measures give answers to questions which should be relevant for a variety of applications.

We have considered the two nonlocal stability measures $\mathcal{P}$ and $\mathcal{D}$,
based on the size and shape of the basin of attraction,
respectively,
These measures account for both large and small perturbations;
$\mathcal{P}$ reflects the probability that the system returns to the attractor given a perturbation,
while $\mathcal{D}$ estimates the smallest perturbation needed to push the solution to another attractor.
The measure $\mathcal{D}$ calls for suitable ways to compare distances in the state-space and, for mechanical systems,
we constructed the version $\mathcal{D}_{\text{energy}}$ reflecting the least amount of energy needed to push the systems solution into another attractor.
We also introduced the version $\mathcal{D}_{\text{rel}}$ considering relative distances when dealing with population dynamics.

We proceeded by considering the systems returntime to the attractor given a perturbation
and defined the simple nonlocal resilience measures $\mathcal{R}$
and $\mathcal{R}_{\text{worst}}$.
In contrary to local resilience measures,
based on linearizations and eigenvalues of the Jacobian matrix,
these measures account for large as well as small perturbations;
$\mathcal{R}$ gives the expected rate of return of the system given a random perturbation while
$\mathcal{R}_{\text{worst}}$ reflects the slowest rate of return of the system given a random perturbation.
Moreover, we introduced the concept basin-time as the subset of the basin of attraction from which all trajectories return to the attractor in time $\tau$.
By replacing the basin of attraction with the smaller set basin-time in the constructions of the stability measures $\mathcal{P}$ and $\mathcal{D}$,
we derived two more resilience measures, $\mathcal{P}^{\tau}$ and $\mathcal{D}^{\tau}$.
The measure $\mathcal{P}^{\tau}$ reflects the probability that the system returns from a random perturbation in time $\tau$,
while $\mathcal{D}^{\tau}$ estimates the least perturbation from which the system will not return in time $\tau$.

We demonstrated our approach on three mathematical models
(the Solow-Swan model of economic growth, 
a system modeling an electro-mechanical machine and 
a stage-structured population model) 
by ``stressing" the dynamics in each model through variations in parameter values.
During these stressed scenarios, we discussed which dynamic behaviour that can be recorded, and why, by each measure.
For example, all nonlocal resilience measures 
detect all considered variations in the Solow-Swan model,
while local measures failed to record the destabilization in three out of four stressed cases.
The measure $\mathcal{D}_{\text{energy}}$ gave the clearest warning signal before the electro-mechanical machine damaged due to a reduction of a spring stiffness, and
all resilience measures agreed on similar results as harvest pressure increased in the stage-structured population model;
they all gave arguments in favour of equal harvesting of juveniles and adults, compared to pure adult harvesting.
In our examples we have seen, and in general we believe,
that the three properties size of basin, shape of basin and returntime constitute a sufficient fundament when building stability and resilience measures.
We can also conclude that no measure is ``the best" measure through all applications 
and thus complementary measures,
such as $\mathcal{P}$ measuring the size of basin and $\mathcal{D}$ measuring the shape of basin,
should be considered simultaneously.

Even though valuable information of a systems stability can easily (using today's computers) be recorded from properties of the basin of attraction and the returntime,
the use of nonlocal measures of stability and resilience are today rare in the literature.
We believe that the simple ideas and measure constructions presented here will constitute an important step in order to quantify stability and resilience through nonlocal methods, and thereby fill this gap in the literature.

%



\end{document}